\documentclass[11pt]{article}

\usepackage{amsmath,amssymb,amsthm}
\usepackage{mathrsfs}
\usepackage{hyperref}
\usepackage{geometry}
\usepackage{array}
\usepackage{booktabs}

\newcommand{\Cl}{\operatorname{Cl}}
\newcommand{\ECl}{\operatorname{ECl}}
\newcommand{\HCl}{\operatorname{HCl}}

\newcommand{\Bcal}{\mathcal{B}}

\newtheorem{proposition}{Proposition}
\newtheorem{corollary}{Corollary}
\newtheorem{remark}{Remark}

\title{Elliptic Clausen Functions and Degenerations\\
\large Circular, Elliptic, and Hyperbolic Parallelism}
\author{Ken Nagai\thanks{Email: \texttt{tknagai@outlook.com}. Independent Researcher.}}
\date{}

\begin{document}
\maketitle

\begin{abstract}
We introduce a unified elliptic extension of CL-type Clausen functions
based on logarithmic primitives of the Jacobi theta function.
The resulting elliptic Clausen family satisfies the same integral recursion
as the classical circular case, with all differences encoded in boundary
constants determined by the underlying logarithmic kernel.
This separation clarifies a strict parallelism between circular, elliptic,
and hyperbolic regimes and makes their degeneration limits transparent.
We further discuss the general structure of the odd boundary constants,
which organize naturally into modular families associated with the elliptic
kernel.
Possible extensions to SL-type frameworks and related master objects are
briefly outlined.
\end{abstract}

\section{Introduction}

Classical Clausen functions arise naturally as primitives of trigonometric
kernel expansions and play a central role in the study of polylogarithms,
Fourier series, and special values of zeta functions.

Classical background on Clausen and polylogarithmic functions
may be found in~\cite{Lewin}.

In this paper, we revisit the \emph{CL-type} Clausen framework and present
a unified elliptic extension whose degeneration limits recover the classical
circular and hyperbolic cases in a fully parallel manner.

Our approach emphasizes integral recursions with explicit boundary constants.
The key observation is that, across circular, elliptic, and hyperbolic regimes,
the recursion structure remains identical, while all differences are absorbed
into boundary terms expressed by logarithmic primitives:
\[
\log\sin,\qquad \log\vartheta_1,\qquad \log\sinh.
\]
This perspective allows a transparent structural comparison and clarifies the
role of degeneration limits.

We deliberately restrict attention to the CL-type framework.
Possible extensions to SL-type Clausen functions, umbral formulations,
and related master objects are deferred to future work.

The paper is organized as follows.
Section~2 recalls the classical CL-type recursion,
Section~3 introduces the elliptic kernel and the corresponding
elliptic Clausen functions,
Sections~4--6 analyze degeneration behavior and the structure
of boundary constants,
and the conclusion summarizes the unified perspective.

\section{Classical Clausen Functions and Integral Recursions}

We briefly recall the classical Clausen functions.
For integers $n\ge1$, they may be defined by
\[
\Cl_n(x) = 
\begin{cases}
\displaystyle \sum_{k=1}^\infty \frac{\sin(kx)}{k^n}, & n \text{ even},\\[1ex]
\displaystyle \sum_{k=1}^\infty \frac{\cos(kx)}{k^n}, & n \text{ odd},
\end{cases}
\]
whenever the series converges.

A fundamental feature of the CL-type family is the integral recursion
\[
\frac{d}{dx}\Cl_{n+1}(x) = \Cl_n(x),
\]
together with boundary constants determined by logarithmic primitives.

Various analytic aspects of classical Clausen functions
have been studied extensively in the literature,
see for example~\cite{Zucker}.

\subsection{Boundary constants in the circular case}

In the circular regime, the primitive of $\Cl_1(x)$ is governed by
\[
\int \Cl_1(x)\,dx = -\log\bigl(2\sin\tfrac{x}{2}\bigr) + \text{const}.
\]
All higher boundary constants inherit this logarithmic structure.

\subsection{Motivation for elliptic extension}

The circular framework closes on trigonometric kernels.
To go beyond this closure while preserving the recursion structure,
it is natural to replace $\log\sin$ by an elliptic logarithmic primitive.
This motivates the introduction of elliptic kernels.

\section{Elliptic Kernel and Elliptic Clausen Functions}

We now pass from the circular setting to an elliptic framework
while preserving the recursion structure described above.

\subsection{Elliptic logarithmic kernel}

Let $\vartheta_1(x\,|\,\tau)$ denote the odd Jacobi theta function.

Standard background on elliptic theta functions
may be found in~\cite{Chandrasekharan}.

We define the elliptic kernel by the logarithmic primitive
\[
K_{\mathrm{ell}}(x;\tau) := \log \vartheta_1(x\,|\,\tau).
\]
This kernel reduces to the circular logarithm under degeneration
$\tau\to i\infty$.

\subsection{Definition of elliptic Clausen functions}

The elliptic Clausen functions $\ECl_n(x;\tau)$ are defined recursively by
\[
\frac{\partial}{\partial x}\ECl_{n+1}(x;\tau) = \ECl_n(x;\tau),
\]
with boundary constants fixed by $K_{\mathrm{ell}}(x;\tau)$.
This definition mirrors the classical CL-type construction,
with the sole modification occurring in the boundary term.

\section{Degenerations and Boundary Constants}

We now describe the degeneration structure.
The recursion itself is invariant; only the boundary constants change.

\begin{proposition}
Across circular, elliptic, and hyperbolic regimes, the integral recursion
for CL-type Clausen functions is identical.
All differences are encoded in the boundary constants, which are given by
\[
\log\sin \quad\leftrightarrow\quad \log\vartheta_1 \quad\leftrightarrow\quad \log\sinh.
\]
\end{proposition}

\begin{remark}
Conceptually, the elliptic regime acts as a master object whose
degenerations reproduce the circular and hyperbolic limits.
\end{remark}

\begin{corollary}
The degeneration limits
\[
\tau\to i\infty, \qquad \tau\to i0^+ \ (\text{via } S\text{-transform})
\]
yield circular and hyperbolic Clausen functions, respectively,
without altering the recursion structure.
\end{corollary}

The modular transformations underlying these degenerations
are classical and discussed in detail in~\cite{Apostol}.

\subsection{Circular degeneration ($\tau\to i\infty$)}

In this limit,
\[
\vartheta_1(x\,|\,\tau) \longrightarrow 2\sin(\pi x),
\]
and $\ECl_n$ reduces to the classical $\Cl_n$.

\subsection{Hyperbolic degeneration ($\tau\to i0^+$)}

After modular transformation,
\[
\vartheta_1(x\,|\,\tau) \longrightarrow \sinh(\pi x),
\]
leading to hyperbolic CL-type functions $\HCl_n$.

This completes the circular–elliptic–hyperbolic parallelism.

\section{Correspondence Table and Structural Overview}

\subsection{Parallel correspondence}

\begin{center}
\begin{tabular}{@{}lll@{}}
\toprule
Regime & Kernel primitive & Boundary constant \\
\midrule
Circular & $\log\sin$ & $\log\sin(\pi x)$ \\
Elliptic & $\log\vartheta_1$ & $\log\vartheta_1(x|\tau)$ \\
Hyperbolic & $\log\sinh$ & $\log\sinh(\pi x)$ \\
\bottomrule
\end{tabular}
\end{center}

\subsection{Structural interpretation}

The table highlights a strict parallelism.
The CL-type recursion defines a common backbone, while each regime
is distinguished solely by its logarithmic kernel.

\medskip
\noindent
\textbf{Summary of boundary constants under degenerations.}
\;
Under the two standard degenerations of the elliptic modulus,
the odd boundary constants $\mathcal B_{2m+1}(\tau)$ reduce to
explicit multiples of the odd zeta values:
\begin{equation}
\label{eq:boundary-degeneration-summary}
\mathcal B_{2m+1}(\tau)\;\longrightarrow\;
\begin{cases}
\zeta(2m+1), & \tau\to i\infty \quad\text{(circular)},\\[6pt]
\dfrac{1}{2^{2m-1}}\zeta(2m+1), & \tau\to i0^+ \quad\text{(hyperbolic)}.
\end{cases}
\end{equation}

\section{Remarks on the General Structure of $\Bcal_{2m+1}(\tau)$}

\begin{remark}
The boundary constants appearing in odd-order recursions naturally assemble
into families $\Bcal_{2m+1}(\tau)$.
Their transformation properties reflect modular behavior inherited from
$\vartheta_1$.
A systematic study of these objects is deferred to future work.
\end{remark}

\section{General structure of the odd boundary constants}
\label{sec:boundary-structure}

In this section we briefly comment on the general structure of the odd boundary
constants $\mathcal B_{2m+1}(\tau)$ appearing in the elliptic Clausen framework.
Although their explicit evaluation is not required for the results of this paper,
their common origin clarifies the conceptual role of boundary terms in the
integral recursions discussed above.

\subsection*{Structural origin}

The boundary constants
\[
\mathcal B_{2m+1}(\tau)
=
\mathrm{ECl}_{2m+1}(0\mid\tau)
\]
arise from a single universal mechanism for all $m$.
They are generated as regularized moments of the elliptic kernel
\begin{equation}
\label{eq:elliptic-kernel}
K_{\mathrm{ell}}(x\mid\tau)
=
-2\log\!\Big(
\frac{\vartheta_1(x\mid\tau)}
{\vartheta_1'(0\mid\tau)\,x}
\Big),
\end{equation}
which naturally appears as the primitive of the elliptic integrand
in the CL-type setting.

Expanding the kernel at the origin,
\begin{equation}
\label{eq:kernel-expansion}
K_{\mathrm{ell}}(x\mid\tau)
=
\sum_{m\ge1} c_{2m}(\tau)\,x^{2m},
\end{equation}
one finds that the constants $\mathcal B_{2m+1}(\tau)$ are determined,
up to universal rational factors, by the coefficients $c_{2m}(\tau)$.
These coefficients play the role of elliptic analogues of the Bernoulli numbers.

From a modular viewpoint, the coefficients $c_{2m}(\tau)$ can be expressed
in terms of classical modular objects, such as Eisenstein series,
or equivalently in terms of Hurwitz-type numbers associated with the
elliptic curve.
In this sense, the family $\{\mathcal B_{2m+1}(\tau)\}_{m\ge0}$ encodes
the modular data of the elliptic kernel in a compact form.

\subsection*{Degenerations}

Under degenerations of the elliptic modulus, this general structure collapses
to familiar constants.
In the circular limit $\tau\to i\infty$, the elliptic kernel
\eqref{eq:elliptic-kernel} reduces to the logarithmic sine kernel,
and the boundary constants $\mathcal B_{2m+1}(\tau)$ converge to
the odd zeta values $\zeta(2m+1)$.
In the hyperbolic degeneration obtained via the modular $S$-transformation
$\tau\mapsto -1/\tau$, the kernel degenerates to the logarithmic hyperbolic
kernel, and the corresponding boundary constants reduce to their
regularized hyperbolic counterparts.

These observations confirm that the difference between the circular,
elliptic, and hyperbolic settings lies entirely in the boundary constants,
while the homogeneous integral recursion itself remains unchanged.

\section{Conclusion and Outlook}

We have presented a unified CL-type framework for Clausen functions
based on elliptic logarithmic kernels.
By isolating boundary constants as the sole source of variation,
we obtained a transparent parallelism between circular, elliptic,
and hyperbolic regimes.

Possible extensions include SL-type Clausen functions,
umbral formulations, and connections to master objects such as $F^\star$.
These directions lie beyond the scope of the present work and
will be addressed elsewhere. This viewpoint suggests that elliptic kernels provide a natural
organizing principle for a broader class of Clausen-type structures.

\section*{Acknowledgments}
The author thanks the collaborative discussions carried out within the \emph{hoge \& fuga} series (2026).

\bibliographystyle{plain}

\begin{thebibliography}{9}

\bibitem{Lewin}
L.~Lewin,
\emph{Polylogarithms and Associated Functions},
North-Holland, 1981.


\bibitem{Apostol}
T.~M.~Apostol,
\emph{Modular Functions and Dirichlet Series in Number Theory},
2nd ed., Springer, 1990.

\bibitem{Zucker}
I.~J.~Zucker,
The evaluation of some definite integrals involving Clausen functions,
\emph{J. Phys. A} \textbf{7} (1974), 1568--1575.

\bibitem{Chandrasekharan}
K.~Chandrasekharan,
\emph{Elliptic Functions},
Springer, 1985.

\end{thebibliography}

\end{document}